\documentclass[11pt,leqno]{article}
\usepackage{amsfonts,amssymb,amsmath,amsthm,latexsym}
\usepackage{color,graphicx}
\usepackage[titletoc]{appendix}

\usepackage{amsfonts} 
\usepackage{latexsym}
\usepackage{graphicx}
\usepackage[utf8]{inputenc}
\usepackage[T1]{fontenc}
\newtheorem{theorem}{\sc Theorem}[section]
\newtheorem{lemma}{\sc Lemma}[section]
\newtheorem{proposition}{\sc Proposition}[section]

\newtheorem{remark}{\sc Remark}[section]

\begin{document}
\title{Unique continuation and time decay for a higher-order water wave model}
\author{ Ademir F. Pazoto \thanks{Institute of Mathematics, Federal University of Rio de Janeiro,
UFRJ,P.O. Box 68530, CEP 21945-970, Rio de Janeiro, RJ, Brasil . E-mail: {\tt ademir@im.ufrj.br}} \and  Miguel D. Soto Vieira
\thanks{Institute of Mathematics, Federal University of Rio de Janeiro, UFRJ, P.O. Box
68530, CEP 21945-970, Rio de Janeiro, RJ, Brasil. E-mail: {\tt dariosv@gmail.com}} }

\date{}

\maketitle

\begin{abstract}
This work is devoted to prove the exponential decay for the energy of solutions of a
higher order Korteweg -de Vries (KdV)--Benjamin-Bona-Mahony (BBM) equation on a periodic domain with a localized damping
mechanism. Following the method in \cite{ro-zhan}, which combines energy estimates, multipliers and compactness arguments,
the problem is reduced to prove the Unique Continuation Property (UCP) for weak solutions of the model. Then, this is done
by deriving Carleman estimates for a system of coupled elliptic-hyperbolic
equations.
\end{abstract}



\section{Introduction}
The field of dispersive equations has received increasing attention since the pioneering works
of Stokes, Boussinesq and Korteweg and de Vries in the nineteenth century. It pertains to a
modern line of research which is important both scientifically and for potential applications.
On the one hand, the mathematical theoretical research of dispersive equations is important
for applied sciences since it has provided solid foundations for the verification and applicability
of these models. On the other hand, this theoretical research has proved to be very valuable
for mathematics itself. Such equations have presented very difficult and interesting challenges,
motivating the development of many new ideas and techniques within mathematical analysis.
Particularly, sets of equations were derived to describe the dynamics of the waves in some specific physical regimes
and much effort has been expended on various aspects of the initial and boundary value problems.
In this context, Bona, Carvajal, Panthee and Scialom \cite{bona-carvajal-panthee-scialom} derived and analyzed
a higher order water wave model to describe the unidirectional propagation of water waves by using the second order approximation
in the two-way model, the so-called $abcd$-system introduced in \cite{BCS1,BCS2}. The model is also known as the fifth order KdV–BBM
type equation and has the form
\begin{equation}\label{Kdv-BBM-1}
\begin{array}{l}
\vspace{2mm}    u_t+u_x-b_1u_{txx}+a_1u_{xxx}+bu_{txxxx}+au_{xxxxx}\\
    \qquad\qquad+\displaystyle\frac{3}{2}uu_x+\gamma(u^2)_{xxx}-\frac{7}{48}(u_{x}^{2})_{x}-\frac{1}{8}(u^3)_x=0.
\end{array}
\end{equation}
The unknown $u$ is a real valued function of the real variables $x$ and $t$ and subscripts indicate partial differentiation. The five parameters $a, a_1, b, b_1$ and $\gamma$ are not arbitrary. Indeed, they are determined by the choice of five more fundamental parameters $\theta, \lambda, \mu, \lambda_1$ and $\mu_1$. The constant $\theta$ has physical significance. It is related to the height at which the horizontal velocity is specified, a dependent variable which does not appear explicitly in these unidirectional models. The parameter $\theta \in [0,1]$ because the vertical coordinate has been scaled by the undisturbed depth $h_0$. The constants $\lambda, \mu, \lambda_1$ and $\mu_1$
are modelling parameters and, in principle, can take any real value. However, the five parameters appearing in \eqref{Kdv-BBM-1} are not independent and should satisfy some relations. All them are shown in \cite{bona-carvajal-panthee-scialom}. Such conditions come from the physics of the problem and we tacitly assume them to hold
throughout the entire paper. Depending on the problem under study (well-posedness of the linear or nonlinear model, stabilization), additional restrictions on the sign of these parameters
will be imposed later on.

The incorporation of damping mechanisms is often crucial in obtaining good agreement between experimental observations and the prediction
of theoretical models describing the propagation of waves in dispersive media. The
problem might be easy to solve when the underlying models have a strong enough intrinsic
dissipative nature, but very often, as the cases we address here, the models are of conservative
nature and the decay requires appropriate damping mechanisms. Obviously, for practical
purposes, it is desirable to achieve this property with a minimal amount of damping both
in what concerns its support and its intensity.

In this work, considerations will be given for model \eqref{Kdv-BBM-1} on the one-dimensional torus $\mathbb{T}=\mathbb{R}/(2\pi\mathbb{Z})$, with a localized damping term. Our purpose is to investigate the dissipative effects generated by this damping when $t\rightarrow \infty$. Our analysis does not depend on particular relations between the coefficients of the system. However, in order to provide the tools needed to deal with the problem, we assume that
		\begin{equation}\label{parameters}
	{\it b, \,\, b_1 > 0.}\nonumber
	\end{equation}
	
Under the above conditions, the following closed-loop system will be considered:
\begin{equation}
\begin{array}{l}\label{nlinitiesta-1}
  u_t+u_x-b_1u_{txx}+a_1u_{xxx}+bu_{txxxx}+au_{xxxxx}\\
\vspace{2mm}  \displaystyle\quad+\frac{3}{2}uu_x+\gamma(u^2)_{xxx}-\frac{7}{48}(u_{x}^{2})_{x}-\frac{1}{8}(u^3)_x =\sigma(I-b_1\partial_{x}^{2}+b\partial_{x}^{4})[\sigma u],\ \ (x,t)\in\mathbb{T}\times(0,T)\\
  u(x,0)=u_0(x), \qquad\qquad\qquad\qquad\qquad\qquad\qquad\qquad\qquad\qquad\qquad\quad x\in\mathbb{T},
\end{array}
\end{equation}
where $\sigma\in \mathcal{C}^\infty(\mathbb{T})$ is a nonzero function. We also introduce the set
\begin{equation}\label{omega}
\omega=\{x\in\mathbb{T}:\sigma(x)\neq 0\}\neq\emptyset.
\end{equation}

Then, if we consider the following equivalent norm in the $L^2$-based Sobolev space $H^2(\mathbb{T})$,
	\begin{equation}\label{energy}
		\|u(\cdot,t)\|_{H^2(\mathbb{T})}^2=\int_{\mathbb{T}}(u^{2} + b_1 u _{x}^{2} + b u_{xx}^{2})dx,
	\end{equation}
at least formally, we can deduce that
	\begin{equation}\label{dtenergy}
		\|u(\cdot,t)\|_{H^2(\mathbb{T})}^2 - \|u_0\|_{H^2(\mathbb{T})}^2=-2\displaystyle\int_0^t\| \sigma u(s)\|_{H^2(\mathbb{T})}^{2}ds.
	\end{equation}
	Indeed, to obtain the identity above we multiply the equation in \eqref{nlinitiesta-1} by $u$ and integrate by parts
	over the spatial domain $\mathbb{T}$. Identity \eqref{dtenergy} shows that the term on the right hand side plays the role of a damping mechanism.
 Moreover, it indicates that the $H^2$-norm decreases along the trajectories of the system, which generates a flow that can be
 continued indefinitely in the temporal variable. Therefore, we can ask whether the solutions converges to zero,
 as $t\rightarrow\infty$, and at which rate they decay.

In view of \eqref{dtenergy}, the problem of the exponential decay can be stated in the following equivalent form: To find $T > 0$
and $C > 0$, such that
\begin{equation}\label{observinequ-1}
    \|u_0\|_{H^2(\mathbb{T})}^{2}\leq C\int_{0}^{T}\|\sigma u(t)\|_{H^2(\mathbb{T})}^{2}dt
\end{equation}
holds for every finite energy solution of \eqref{nlinitiesta-1}. Indeed, from \eqref{observinequ-1} and \eqref{dtenergy}, we have that $E(T ) \leq \beta E(0)$ with
$0 < \beta < 1$, which combined with the semigroup property, allow us to derive the exponential decay of the $H^2$-norm of the solutions.

This paper is devoted to analyze this problem.
Our analysis was inspired by the results proved in \cite{ro-zhan} from which we borrow the main ideas involved in our proofs.
In fact, proceeding as in \cite{ro-zhan}, which combine multipliers, energy estimates and compactness arguments, the problem of obtaining \eqref{observinequ-1}
is reduced to show that the unique solution of \eqref{nlinitiesta-1}, such that $\sigma(x)u = 0$ everywhere, has to be the trivial one.
This problem may be viewed as a unique continuation
one since $ \sigma u = 0$ implies that $u = 0$ in $\omega\times (0,T)$, for $\omega$ given by \eqref{omega}.
To solve this problem we develop a Carleman inequalities for weak solutions of the model as it has been done
in \cite{ro-zhan} for a KdV-BBM equation. The equation is first split into a coupled system of an elliptic equation and a transport equation.
Next, some Carleman estimates are derived with the same singular weights for both the elliptic and the hyperbolic equations. Finally, the unique continuation
result is proved by combining these Carleman estimates with a regularization process. However, this unique continuation
result cannot be applied directly due to the regularity of solutions we are dealing with. Therefore, some additional assumptions concerning
the initial conditions are needed. More precisely, we assume that $0<r\leq \|u_0\|_{H^2(\mathbb{T})}$ and $\|u_0\|_{H^3(\mathbb{T})}\leq R$, for given $r, R >0$.
These conditions seem to be technical, but, in their absence, we do not know how to derive the result. It is also important to emphasize that the choice of the damping
function was motivated by the controllability results obtained in \cite{bau-pazo} for corresponding linearized system. Indeed, as an application of the
controllability result, we construct feedback controls via some results obtained in \cite{LIU,SLEM}, such that the resulting linearized
closed-loop system is shown to be exponentially stable. Such stabilization results plays a crucial role in the analysis we describe above.

We remark that the unique continuation property for the BBM equation is still an open problem. Moreover, since the underlying Cauchy problem
is a characteristic one, we can not expect to apply Carleman-type estimates or the classical Holmgren uniqueness theorem. In \eqref{Kdv-BBM-1},
the presence of the higher-order KdV term $u_{xxxxx}$ results in  much better properties and allows to establish a unique continuation result.

It is also worth mentioning that, in addition to the analysis carried out in the work, the stabilization of the KdV and BBM
equations have also been studied considering other types of damping mechanisms. In this sense, we refer to \cite{rosier-zhang surv,ro-zhan}
for a quite complete review of the field.

The remainder of this paper is organized as follows: in Section 2, we prove some Carleman estimates and derive a unique continuation property for the higher-order
KdV-BBM equation. Section 3 is devoted to the stabilization of the damped equation. Finally, for the sake of completeness, we include in an Appendix
some computations used in Section 3 and the proof a Carleman estimate derived in \cite{ro-zhan}.

\section{Carleman estimates}
The first part of this section is devoted to prove an appropriate Carleman estimate for the higher-order KdV-BBM equation
\begin{equation}
\begin{array}{l}\label{Kdv-BBM}
\vspace{2mm}  u_t+u_x-b_1u_{txx}+a_1u_{xxx}+bu_{txxxx}+au_{xxxxx}\\
\displaystyle\qquad\qquad+\frac{3}{2}uu_x+\gamma(u^2)_{xxx}-\frac{7}{48}(u_{x}^{2})_{x}-\frac{1}{8}(u^3)_x =0,\quad (x,t)\in\mathbb{T}\times(0,T),\\
\end{array}
\end{equation}
or
\begin{equation}\label{carleq}
     u_t-b_1u_{txx}+bu_{txxxx}+au_{xxxxx}+q(u)u_x+p(u)u_{xxx}+r(u)u_{xx}=0, \ \ \ (x,t)\in\mathbb{T}\times(0,T),
\end{equation}
where $q(u)=1+\frac{3}{2}u-\frac{3}{8}u^2$, $p(u)=a_1+2\gamma u$ and $r(u)=(6\gamma-\frac{7}{24})u_x$.

\vglue 0.2 cm

Next, this Carleman estimate is employed to prove the following unique continuation result:
\begin{theorem}\label{ucp}
Let $a,b\neq0$, $T>\frac{2\pi b}{|a|}$ and $q,p,r\in L^{\infty}(0,T;L^{\infty}(\mathbb{T}))$. Let $\omega\subset\mathbb{T}$ be a nonempty open set. Let $u\in L^2(0,T;H^4(\mathbb{T}))\cup L^{\infty}(0,T;H^3(\mathbb{T}))$ satisying \eqref{carleq} and
\begin{align}\label{contuni}
    u(x,t)=0 \  \ \ \text{for a.e}\ \ (x,t)\in \omega\times(0,T).
\end{align}
Then, $u\equiv 0$ in $\mathbb{T}\times(0,T)$.
\end{theorem}
\begin{proof}
Assume that
\begin{align}
u\in L^2(0,T;H^4(\mathbb{T}))
\end{align}
and let $w=u-b_1u_{xx}+bu_{xxxx}\in L^2(0,T;L^2(\mathbb{T}))$. Then, the pair $(u,w)$ solves the following system
\begin{align}
    &u-b_1u_{xx}+bu_{xxxx}=w \label{elitica}\\
    &w_t+\frac{a}{b}w_x=(\frac{a}{b}-q)u_x-(\frac{ab_1}{b}+p)u_{xxx}-ru_{xx}.\label{transpor}
\end{align}
In order to make more clear the assumption of the theorem, we recall the following remark, which can be found in \cite{ro-zhan}:
\begin{remark}
There is a finite speed propagation for KdV-BBM. For instance, if we assume that $q(x)=\frac{a}{b}$, $p(x)=\frac{-ab_1}{b}$ and $r(x)=0$ for all $x\in\mathbb{T}$, where $a,b,b_1>0$ are given, and that $\omega=(2\pi-\epsilon,2\pi)$ for a small $\epsilon>0$, then the UCP fails in time $T\leq \frac{b(2\pi-2\epsilon)}{a}$. Indeed, picking any nontrivial initial state $u_0\in C^{\infty}_{0}(0,\epsilon)$, we easily see that the solution $(u,w)$ of \eqref{elitica}-\eqref{transpor} is $u(x,t)=u_0(x-\frac{a}{b}t)$, $w(x,t)=w_0(x-\frac{a}{b}t)$, where $w_0=(I-\partial_{x}^{2})u_0.$ Then, $u(x,t)=0$  for $(x,t)\in\omega\times(0,\frac{b(2\pi-2\epsilon)}{a})$ although $u\neq 0$. Hence, the condition $T>\frac{2b\pi}{|a|}$ in the Theorem \ref{ucp} is sharp.
\end{remark}

We choose a system of coordinates and indentify $\mathbb{T}$ with $[0,2\pi)$. Thus, without loss of generality, we can assume that $a>0$, and that $\omega=(2\pi-\eta,2\pi+\eta)\sim [0,\eta)\cup (2\pi-\eta,2\pi)$ for some $\eta\in(0,\pi)$. Then, we consider $T$, such that
\begin{align}
    T>\frac{2b\pi}{a}.
\end{align}

Following the approach developed in \cite{ro-zhan}, we first obtain some Carleman estimates for the elliptic equation \eqref{elitica} and the transport equation \eqref{transpor} with the same weights function. Next, we combine them to derive a single one for \eqref{carleq}. In order to do that, we pick $\delta>0$ and $\rho\in(0,1)$, such that
\begin{align}\label{deltaro}
    \frac{\rho a}{b} T>2\pi+\delta
\end{align}
and a function $\psi\in C^{\infty}([0,2\pi]\times \mathbb{R})$ satisfying
\begin{align}
    &\psi(x)=(x+\delta)^2\ \ \ \ \ \ \ \ \ \ \ \ \ \text{for}\ x\in[\frac{\eta}{2},2\pi-\frac{\eta}{2}],\label{poli}\\
    &\frac{d^k\psi}{dx^k}(0)=\frac{d^k\psi}{dx^k}(2\pi)\ \ \ \ \ \ \ \ \text{for}\ k=1,2,3,4,5,6,7,\label{condperi}\\
    &2\delta\leq\frac{d\psi}{dx}(x)\leq2(2\pi+\delta)\ \ \text{for}\ x\in[0,2\pi].\label{limitderiv}
\end{align}
Then, we introduce the function $\varphi\in C^{\infty}([0,2\pi]\times\mathbb{R})$ as follows
\begin{align}\label{varphi}
    \varphi(x,t)=\psi(x)-\rho a^2t^2.
\end{align}

Under the above considerations, we derive the following Carleman estimate for \eqref{carleq}.
\begin{proposition}\label{carleman}
Let $\omega$, $a$ and $T$ be as above. Then, there exists some positive numbers $s_2$ and $C_2$, such that, for all $s\geq s_2$ and all $u\in L^2(0,T;H^4(\mathbb{T}))$ satisfying \eqref{carleq}, we have
\begin{align}\label{carlecomple}
    \int_{0}^{T}\int_{\mathbb{T}}[s|u_{xxxx}|^2+&s|u_{xxx}|^2+s^3|u_{xx}|^2+s^5|u_x|^2+s^7|u|^2]e^{2s\varphi}dxdt\nonumber\\
     &+s\int_{\mathbb{T}}[|u-b_1u_{xx}+bu_{xxxx}|^2e^{2s\varphi}]_{t=0}dx\\
    &\leq C_2\int_{0}^{T}\int_{\omega}[s|u_{xxxx}|^2+s^3|u_{xx}|^2+s^7|u|^2]e^{2s\varphi}dxdt.\nonumber
\end{align}
\end{proposition}
\begin{proof} In order to make the reading easier, we proceed in several steps.

\vglue 0.2 cm

\noindent $\bullet$ Carleman estimate for the transport equation.

\begin{lemma}\label{carltra}
There exist $s_1\geq s_0$  and $C_1>0$ such that for all $s\geq s_1$ and all $w\in L^2(\mathbb{T}\times(0,T))$ with $w_t+\frac{a}{b}w_x\in L^2(\mathbb{T}\times(0,T))$, the following holds
\begin{equation}\label{Carlemantran}
\begin{array}{l}
    \displaystyle\int_{0}^{T}\int_{\mathbb{T}}s|w|^2e^{2s\varphi}dxdt+\int_{\mathbb{T}}s[|w|^2e^{2s\varphi}]_{t=0}dx+\int_{\mathbb{T}}s[|w|^2e^{2s\varphi}]_{t=T}dx\\
    \leq C_1\displaystyle\left(\int_{0}^{T}\int_{\mathbb{T}}|w_t+\frac{a}{b}w_x|^2e^{2s\varphi}dxdt+\int_{0}^{T}\int_{\omega}s|w|^2e^{2s\varphi}dxdt\right).
\end{array}
\end{equation}
\end{lemma}
\begin{proof}
The result was proved in \cite[Lemma 5.5]{ro-zhan}. For the sake of completeness we have included the proof in the Appendix.
\end{proof}

\noindent $\bullet$ Carleman estimate for the elliptic equation.
\begin{lemma}\label{Carleli}
There exist $s_0\geq1$ and $C_0>0$, such that, for all $s\geq s_0$ and all $u\in H^4(\mathbb{T})$, the following holds
\begin{equation}\label{carlemanelitica}
\begin{array}{l}
    \displaystyle\int_{\mathbb{T}}[s|u_{xxx}|^2+s^3|u_{xx}|^2+s^5|u_x|^2+s^7|u|^2]e^{2s\psi}dx\\
    \qquad\leq C_0\displaystyle\left(\int_{\mathbb{T}}|u_{xxxx}|^2e^{2s\psi}dx+\int_{\omega}(s^7|u|^2+s^3|u_{xx}|^2)e^{2s\psi}dx\right).
    \end{array}
\end{equation}
\end{lemma}
\begin{proof}
We start by considering a classical change of function and a appropriate differential operator. More precisely, let $v=e^{s\psi}u$ and $P=\partial_{x}^{4}$. Then,
we decompose $e^{s\psi}Pu$ as follows
\begin{align*}
    e^{s\psi}Pu=e^{s\psi}P(e^{-s\psi}v)=P_pv+P_nv,
\end{align*}
where
\begin{align}
    &P_pv=(s^4\psi_{x}^{4}+3s^2\psi_{xx}^{2}+4s^2\psi_{xxx}\psi_x)v+12s^2\psi_x\psi_{xx}v_x+6s^2\psi_{x}^{2}v_{xx}+v_{xxxx},\label{partepositiva}\\
    &P_nv=-(6s^3\psi_{x}^{2}\psi_{xx}+s\psi_{xxxx})v-(4s^3\psi_{x}^{3}+4s\psi_{xxx})v_x-6s\psi_{xx}v_{xx}-4s\psi_xv_{xxx}.\label{partenegativa}
\end{align}
From the above decomposition, we obtain
\begin{align*}
    \|e^{s\psi}Pu\|^2=\|P_pv\|^2+\|P_nv\|^2+2(P_pv,P_nv),
\end{align*}
where $(f,g)=\int_\mathbb{T}fgdx$, and $\|f\|=(f,f)$.

The next steps are devoted to analyze the inner product in the identity above. First, observe that,
\eqref{partepositiva} and \eqref{partenegativa} lead to
\begin{align*}
    &(P_pv,P_nv)=((s^4\psi_{x}^{4}+3s^2\psi_{xx}^{2}+4s^2\psi_{xxx}\psi_x)v,-(6s^3\psi_{x}^{2}\psi_{xx}+s\psi_{xxxx})v)\\
    &+((s^4\psi_{x}^{4}+3s^2\psi_{xx}^{2}+4s^2\psi_{xxx}\psi_x)v,-(4s^3\psi_{x}^{3}+4s\psi_{xxx})v_x)\\
    &+((s^4\psi_{x}^{4}+3s^2\psi_{xx}^{2}+4s^2\psi_{xxx}\psi_x)v,-6s\psi_{xx}v_{xx})\\
    &+((s^4\psi_{x}^{4}+3s^2\psi_{xx}^{2}+4s^2\psi_{xxx}\psi_x)v,-4s\psi_xv_{xxx})+(12s^2\psi_x\psi_{xx}v_x,-(6s^3\psi_{x}^{2}\psi_{xx}+s\psi_{xxxx})v)\\
    &+(12s^2\psi_x\psi_{xx}v_x,-(4s^3\psi_{x}^{3}+4s\psi_{xxx})v_x)+(12s^2\psi_x\psi_{xx}v_x,-6s\psi_{xx}v_{xx})\\
    &+(12s^2\psi_x\psi_{xx}v_x,-4s\psi_xv_{xxx})+(6s^2\psi_x^{2}v_{xx},-(6s^3\psi_{x}^{2}\psi_{xx}+s\psi_{xxxx})v)\\
    &+(6s^2\psi_x^{2}v_{xx},-(4s^3\psi_{x}^{3}+4s\psi_{xxx})v_x)+(6s^2\psi_x^{2}v_{xx},-6s\psi_{xx}v_{xx})+(6s^2\psi_x^{2}v_{xx},-4s\psi_xv_{xxx})\\
    &+(v_{xxxx},-(6s^3\psi_{x}^{2}\psi_{xx}+s\psi_{xxxx})v)+(v_{xxxx},-(4s^3\psi_{x}^{3}+4s\psi_{xxx})v_x)+(v_{xxxx},-6s\psi_{xx}v_{xx})\\
    &+(v_{xxxx},-4s\psi_xv_{xxx})\\
    &=\sum_{n=1}^{16}I_n.
\end{align*}
Next, we compute each term $I_n$. This is done employing \eqref{condperi} and integration by parts in $x$:
\begin{align*}
    I_1&=-\int_{\mathbb{T}}(s^4\psi_{x}^{4}+3s^2\psi_{xx}^{2}+4s^2\psi_{xxx}\psi_x)(6s^3\psi_{x}^{2}\psi_{xx}+s\psi_{xxxx})v^2dx\\
    I_2&=-\int_{\mathbb{T}}(s^4\psi_{x}^{4}+3s^2\psi_{xx}^{2}+4s^2\psi_{xxx}\psi_x)(4s^3\psi_{x}^{3}+4s\psi_{xxx})v_xvdx\\
    &=\int_{\mathbb{T}}[(s^4\psi_{x}^{4}+3s^2\psi_{xx}^{2}+4s^2\psi_{xxx}\psi_x)(2s^3\psi_{x}^{3}+2s\psi_{xxx})]_xv^{2}dx\\
    I_3&=-\int_{\mathbb{T}}(s^4\psi_{x}^{4}+3s^2\psi_{xx}^{2}+4s^2\psi_{xxx}\psi_x)(6s\psi_{xx})vv_{xx}dx\\
    &=\int_{\mathbb{T}}[(s^4\psi_{x}^{4}+3s^2\psi_{xx}^{2}+4s^2\psi_{xxx}\psi_x)(6s\psi_{xx})]_xvv_xdx\\
    &+\int_{\mathbb{T}}[(s^4\psi_{x}^{4}+3s^2\psi_{xx}^{2}+4s^2\psi_{xxx}\psi_x)(6s\psi_{xx})]v_{x}^{2}dx\\
    &=-\int_{\mathbb{T}}[(s^4\psi_{x}^{4}+3s^2\psi_{xx}^{2}+4s^2\psi_{xxx}\psi_x)(3s\psi_{xx})]_{xx}v^2dx\\
    &+\int_{\mathbb{T}}[(s^4\psi_{x}^{4}+3s^2\psi_{xx}^{2}+4s^2\psi_{xxx}\psi_x)(6s\psi_{xx})]v_{x}^{2}dx
    \end{align*}
    \begin{align*}
    I_4&=-\int_{\mathbb{T}}(s^4\psi_{x}^{4}+3s^2\psi_{xx}^{2}+4s^2\psi_{xxx}\psi_x)4s\psi_xvv_{xxx}dx\\
    &=\int_{\mathbb{T}}[(s^4\psi_{x}^{4}+3s^2\psi_{xx}^{2}+4s^2\psi_{xxx}\psi_x)4s\psi_x]_xvv_{xx}dx\\
    &+\int_{\mathbb{T}}(s^4\psi_{x}^{4}+3s^2\psi_{xx}^{2}+4s^2\psi_{xxx}\psi_x)4s\psi_xv_xv_{xx}dx\\
    &=-\int_{\mathbb{T}}[(s^4\psi_{x}^{4}+3s^2\psi_{xx}^{2}+4s^2\psi_{xxx}\psi_x)4s\psi_x]_{xx}vv_{x}dx\\
    &-\int_{\mathbb{T}}[(s^4\psi_{x}^{4}+3s^2\psi_{xx}^{2}+4s^2\psi_{xxx}\psi_x)4s\psi_x]_xv_{x}^{2}dx\\
    &-\int_{\mathbb{T}}[(s^4\psi_{x}^{4}+3s^2\psi_{xx}^{2}+4s^2\psi_{xxx}\psi_x)2s\psi_x]_xv_{x}^{2}dx\\
    &=\int_{\mathbb{T}}[(s^4\psi_{x}^{4}+3s^2\psi_{xx}^{2}+4s^2\psi_{xxx}\psi_x)2s\psi_x]_{xxx}v^2dx\\
    &-\int_{\mathbb{T}}[(s^4\psi_{x}^{4}+3s^2\psi_{xx}^{2}+4s^2\psi_{xxx}\psi_x)6s\psi_x]_xv_{x}^{2}dx\\
    I_5&=-\int_{\mathbb{T}}12s^2\psi_x\psi_{xx}(6s^3\psi_{x}^{2}\psi_{xx}+s\psi_{xxxx})vv_xdx\\
    &=\int_{\mathbb{T}}[6s^2\psi_x\psi_{xx}(6s^3\psi_{x}^{2}\psi_{xx}+s\psi_{xxxx})]_xv^2dx\\
    I_6&=-\int_{\mathbb{T}}12s^2\psi_x\psi_{xx}(4s^3\psi_{x}^{3}+4s\psi_{xxx})v_{x}^{2}dx\\
    I_7&=-\int_{\mathbb{T}}12s^2\psi_x\psi_{xx}6s\psi_{xx}v_xv_{xx}dx=\int_{\mathbb{T}}[36s^3\psi_x\psi_{xx}^{2}]_xv_{x}^{2}dx\\
  I_8&=-\int_{\mathbb{T}}12s^2\psi_x\psi_{xx}4s\psi_{x}v_xv_{xxx}dx=-\int_{\mathbb{T}}48s^3\psi_{x}^{2}\psi_{xx}v_{x}v_{xxx}dx\\
  &=\int_{\mathbb{T}}[48s^3\psi_{x}^{2}\psi_{xx}]_xv_{x}v_{xx}dx+\int_{\mathbb{T}}48s^3\psi_{x}^{2}\psi_{xx}v_{xx}^{2}dx\\
  &=-\int_{\mathbb{T}}[24s^3\psi_{x}^{2}\psi_{xx}]_{xx}v_{x}^{2}dx+\int_{\mathbb{T}}48s^3\psi_{x}^{2}\psi_{xx}v_{xx}^{2}dx\\
  I_9&=-\int_{\mathbb{T}}6s^2\psi_x^{2}(6s^3\psi_{x}^{2}\psi_{xx}+s\psi_{xxxx})vv_{xx}dx\\
  &=\int_{\mathbb{T}}[6s^2\psi_x^{2}(6s^3\psi_{x}^{2}\psi_{xx}+s\psi_{xxxx})]_xvv_{x}dx+\int_{\mathbb{T}}6s^2\psi_x^{2}(6s^3\psi_{x}^{2}\psi_{xx}+s\psi_{xxxx})v_{x}^{2}dx\\
  &=-\int_{\mathbb{T}}[3s^2\psi_x^{2}(6s^3\psi_{x}^{2}\psi_{xx}+s\psi_{xxxx})]_{xx}v^2dx+\int_{\mathbb{T}}6s^2\psi_x^{2}(6s^3\psi_{x}^{2}\psi_{xx}+s\psi_{xxxx})v_{x}^{2}dx\\
  I_{10}&=-\int_{\mathbb{T}}6s^2\psi_x^{2}(4s^3\psi_{x}^{3}+4s\psi_{xxx})v_{x}v_{xx}dx=\int_{\mathbb{T}}[3s^2\psi_x^{2}(4s^3\psi_{x}^{3}+4s\psi_{xxx})]_xv_{x}^{2}dx\\
  I_{11}&=-\int_{\mathbb{T}}6s^2\psi_x^{2}6s\psi_{xx}v_{xx}^{2}dx=-\int_{\mathbb{T}}36s^3\psi_x^{2}\psi_{xx}v_{xx}^{2}dx\\
  I_{12}&=-\int_{\mathbb{T}}6s^2\psi_x^{2}4s\psi_{x}v_{xx}v_{xxx}dx=\int_{\mathbb{T}}[12s^3\psi_x^{3}]_xv_{xx}^{2}\\
    I_{13}&=-\int_{\mathbb{T}}(6s^3\psi_{x}^{2}\psi_{xx}+s\psi_{xxxx})vv_{xxxx}dx=\int_{\mathbb{T}}[6s^3\psi_{x}^{2}\psi_{xx}+s\psi_{xxxx})]_xvv_{xxx}dx\\
  &+\int_{\mathbb{T}}(6s^3\psi_{x}^{2}\psi_{xx}+s\psi_{xxxx})v_{x}v_{xxx}dx
  \end{align*}
    \begin{align*}
  &=-\int_{\mathbb{T}}[6s^3\psi_{x}^{2}\psi_{xx}+s\psi_{xxxx})]_{xx}vv_{xx}dx-\int_{\mathbb{T}}[6s^3\psi_{x}^{2}\psi_{xx}+s\psi_{xxxx})]_xv_xv_{xx}dx\\
    &-\int_{\mathbb{T}}[6s^3\psi_{x}^{2}\psi_{xx}+s\psi_{xxxx})]_xv_{x}v_{xx}dx-\int_{\mathbb{T}}(6s^3\psi_{x}^{2}\psi_{xx}+s\psi_{xxxx})v_{xx}^{2}dx\\
  &=\int_{\mathbb{T}}[6s^3\psi_{x}^{2}\psi_{xx}+s\psi_{xxxx})]_{xxx}vv_{x}dx+\int_{\mathbb{T}}[6s^3\psi_{x}^{2}\psi_{xx}+s\psi_{xxxx})]_{xx}v_{x}^{2}dx\\
  &+\int_{\mathbb{T}}[6s^3\psi_{x}^{2}\psi_{xx}+s\psi_{xxxx})]_{xx}v_{x}^{2}dx-\int_{\mathbb{T}}(6s^3\psi_{x}^{2}\psi_{xx}+s\psi_{xxxx})v_{xx}^{2}dx\\
  &=-\int_{\mathbb{T}}\frac{1}{2}[6s^3\psi_{x}^{2}\psi_{xx}+s\psi_{xxxx})]_{xxx}v^2dx+\int_{\mathbb{T}}2[6s^3\psi_{x}^{2}\psi_{xx}+s\psi_{xxxx})]_{xx}v_{x}^{2}dx\\
  &-\int_{\mathbb{T}}(6s^3\psi_{x}^{2}\psi_{xx}+s\psi_{xxxx})v_{xx}^{2}dx\\
  I_{14}&=-\int_{\mathbb{T}}(4s^3\psi_{x}^{3}+4s\psi_{xxx})v_xv_{xxxx}dx=\int_{\mathbb{T}}[(4s^3\psi_{x}^{3}+4s\psi_{xxx})]_xv_xv_{xxx}dx\\
  &+\int_{\mathbb{T}}(4s^3\psi_{x}^{3}+4s\psi_{xxx})v_{xx}v_{xxx}dx\\
  &=-\int_{\mathbb{T}}[(4s^3\psi_{x}^{3}+4s\psi_{xxx})]_{xx}v_xv_{xx}dx-\int_{\mathbb{T}}[(4s^3\psi_{x}^{3}+4s\psi_{xxx})]_xv_{xx}^{2}dx\\
  &-\int_{\mathbb{T}}[(2s^3\psi_{x}^{3}+2s\psi_{xxx})]_xv_{xx}^{2}dx\\
  &=\int_{\mathbb{T}}[(2s^3\psi_{x}^{3}+2s\psi_{xxx})]_{xxx}v_{x}^{2}dx-\int_{\mathbb{T}}3[(2s^3\psi_{x}^{3}+2s\psi_{xxx})]_xv_{xx}^{2}dx\\
  I_{15}&=-\int_{\mathbb{T}}6s\psi_{xx}v_{xx}v_{xxxx}dx=\int_{\mathbb{T}}6s\psi_{xxx}v_{xx}v_{xxx}dx+\int_{\mathbb{T}}6s\psi_{xx}v_{xxx}^{2}dx\\
  &=-\int_{\mathbb{T}}3s\psi_{xxxx}v_{xx}^{2}dx+\int_{\mathbb{T}}6s\psi_{xx}v_{xxx}^{2}dx\\
  I_{16}&=-\int_{\mathbb{T}}4s\psi_xv_{xxx}v_{xxxx}dx=\int_{\mathbb{T}}2s\psi_{xx}v_{xxx}^{2}dx.
\end{align*}
Combining the identities above, we get
\begin{align*}
    \|e^{s\psi}Pu\|^2&=\|P_pv\|^2+\|P_nv\|^2+2\int_{\mathbb{T}}h_1(\psi)v^2dx+2\int_{\mathbb{T}}h_2(\psi)v_{x}^{2}dx+\\
    &2\int_{\mathbb{T}}h_3(\psi)v_{xx}^{2}dx+2\int_{\mathbb{T}}h_4(\psi)v_{xxx}^{2}dx,
\end{align*}
where
\begin{align*}
    h_1(\psi)&=[(s^4\psi_{x}^{4}+3s^2\psi_{xx}^{2}+4s^2\psi_{xxx}\psi_x)(2s^3\psi_{x}^{3}+2s\psi_{xxx})]_x\\
    &-(s^4\psi_{x}^{4}+3s^2\psi_{xx}^{2}+4s^2\psi_{xxx}\psi_x)(6s^3\psi_{x}^{2}\psi_{xx}+s\psi_{xxxx})\\
    &-[(s^4\psi_{x}^{4}+3s^2\psi_{xx}^{2}+4s^2\psi_{xxx}\psi_x)(3s\psi_{xx})]_{xx}\\
    &+[(s^4\psi_{x}^{4}+3s^2\psi_{xx}^{2}+4s^2\psi_{xxx}\psi_x)2s\psi_x]_{xxx}\\
    &+[12s^2\psi_x\psi_{xx}(6s^3\psi_{x}^{2}\psi_{xx}+s\psi_{xxxx})]_x-[3s^2\psi_x^{2}(6s^3\psi_{x}^{2}\psi_{xx}+s\psi_{xxxx})]_{xx}\\
    &-\frac{1}{2}[6s^3\psi_{x}^{2}\psi_{xx}+s\psi_{xxxx})]_{xxx}
    \end{align*}
    \begin{align*}
    h_2(\psi)&=[(s^4\psi_{x}^{4}+3s^2\psi_{xx}^{2}+4s^2\psi_{xxx}\psi_x)(6s\psi_{xx})]-[(s^4\psi_{x}^{4}+3s^2\psi_{xx}^{2}+4s^2\psi_{xxx}\psi_x)6s\psi_x]_x\\
    &-12s^2\psi_x\psi_{xx}(4s^3\psi_{x}^{3}+4s\psi_{xxx})+[36s^3\psi_x\psi_{xx}^{2}]_x-[24s^3\psi_{x}^{2}\psi_{xx}]_{xx}\\
    &+[3s^2\psi_x^{2}(4s^3\psi_{x}^{3}+4s\psi_{xxx})]_x\\
    &+2[6s^3\psi_{x}^{2}\psi_{xx}+s\psi_{xxxx})]_{xx}+[(2s^3\psi_{x}^{3}+2s\psi_{xxx})]_{xxx}+6s^2\psi_x^{2}(6s^3\psi_{x}^{2}\psi_{xx}+s\psi_{xxxx})\\
    h_3(\psi)&=48s^3\psi_{x}^{2}\psi_{xx}-36s^3\psi_x^{2}\psi_{xx}+[12s^3\psi_x^{3}]_x-(6s^3\psi_{x}^{2}\psi_{xx}+s\psi_{xxxx})\\
    &-3[(2s^3\psi_{x}^{3}+2s\psi_{xxx})]_x-3s\psi_{xxxx}\\
    h_4(\psi)&=6s\psi_{xx}+2s\psi_{xx}=8s\psi_{xx}.
\end{align*}
The choice of the function $\psi$ given by \eqref{poli} allows us to conclude that there exist $s_0\geq1$, $K>0$ and $K_1>0$, such that, for all $s\geq s_0$,
\begin{align*}
    2h_1(\psi)\geq Ks^7 \ \ \ \ \text{for}\ \ \  (x,t)\in(\frac{\eta}{2},2\pi-\frac{\eta}{2})\times(0,T),\\
    2h_2(\psi)\geq Ks^5 \ \ \ \ \text{for}\ \ \  (x,t)\in(\frac{\eta}{2},2\pi-\frac{\eta}{2})\times(0,T),\\
    2h_3(\psi)\geq Ks^3 \ \ \ \ \text{for}\ \ \  (x,t)\in(\frac{\eta}{2},2\pi-\frac{\eta}{2})\times(0,T), \\
    2h_4(\psi)\geq Ks \ \ \ \ \text{for}\ \ \  (x,t)\in(\frac{\eta}{2},2\pi-\frac{\eta}{2})\times(0,T),
    \end{align*}
and, if $x\in \omega_0=[0,\frac{\eta}{2})\cup(2\pi-\frac{\eta}{2},2\pi)$,
\begin{align*}
    |2h_1(\psi)|\leq K_1s^7\ \ \ \ \text{for}\ \ \  (x,t)\in \omega_0\times(0,T),\\
    |2h_2(\psi)|\leq K_1s^5\ \ \ \ \text{for}\ \ \  (x,t)\in \omega_0\times(0,T),\\
    |2h_3(\psi)|\leq K_1s^3\ \ \ \ \text{for}\ \ \  (x,t)\in \omega_0\times(0,T),\\
    |2h_4(\psi)|\leq K_1s\ \ \ \ \text{for}\ \ \  (x,t)\in \omega_0\times(0,T).
\end{align*}
Consequently, for $s\geq s_0$, we obtain $C>0$ satisfying
\begin{align*}
&\|P_pv\|^2+\int_{\mathbb{T}}[s^7|v|^2+s^5|v_x|^2+s^3|v_{xx}|^2+s|v_{xxx}|^2]dx\\
&=\|P_pv\|^2+\int_{\omega_0}[s^7|v|^2+s^5|v_x|^2+s^3|v_{xx}|^2+s|v_{xxx}|^2]dx\\
&+\int_{\mathbb{T}\backslash\omega_0}[s^7|v|^2+s^5|v_x|^2+s^3|v_{xx}|^2+s|v_{xxx}|^2]dx\\
&\leq\|P_pv\|^2+\int_{\omega_0}[s^7|v|^2+s^5|v_x|^2+s^3|v_{xx}|^2+s|v_{xxx}|^2]dx\\  &+C\int_{\mathbb{T}\backslash\omega_0}[2h_1(\psi)|v|^2+2h_2(\psi)|v_x|^2+2h_3(\psi)|v_{xx}|^2+2h_4(\psi)|v_{xxx}|^2]dx\\
&\leq \|P_pv\|^2+\int_{\omega_0}[s^7|v|^2+s^5|v_x|^2+s^3|v_{xx}|^2+s|v_{xxx}|^2]dx\\  &+C\int_{\mathbb{T}}[2h_1(\psi)|v|^2+2h_2(\psi)|v_x|^2+2h_3(\psi)|v_{xx}|^2+2h_4(\psi)|v_{xxx}|^2]dx\\
&\leq C\left(\|e^{s\psi}Pu\|^2+\int_{\omega_0}[s^7|v|^2+s^5|v_x|^2+s^3|v_{xx}|^2+s|v_{xxx}|^2]dx\right),
\end{align*}
that is,
\begin{equation}\label{precarleelliti}
\begin{array}{l}
  \|P_pv\|^2+\displaystyle\int_{\mathbb{T}}[s^7|v|^2+s^5|v_x|^2+s^3|v_{xx}|^2+s|v_{xxx}|^2]dx\\
  \qquad\leq C\displaystyle\left(\|e^{s\psi}Pu\|^2+\int_{\omega_0}[s^7|v|^2+s^5|v_x|^2+s^3|v_{xx}|^2+s|v_{xxx}|^2]dx\right).
\end{array}
\end{equation}
Observe that $\int_{\mathbb{T}}s^{-1}|v_{xxxx}|^2dx$ is less than the  hand side of \eqref{precarleelliti}. Indeed,
\begin{align*}
    &\int_{\mathbb{T}}s^{-1}|v_{xxxx}|^2dx=\int_{\mathbb{T}}s^{-1}|P_pv-(s^4\psi_{x}^{4}+3s^2\psi_{xx}^{2}+4s^2\psi_{xxx}\psi_x)v-12s^2\psi_x\psi_{xx}v_x-6s^2\psi_{x}^{2}v_{xx}|^2dx\\
    &\leq C\int_{\mathbb{T}}s^{-1}(|P_pv|^2+|s^4\psi_{x}^{4}+3s^2\psi_{xx}^{2}+4s^2\psi_{xxx}\psi_x|^2|v|^2+|12s^2\psi_x\psi_{xx}|^2|v_x|^2+|6s^2\psi_{x}^{2}|^2|v_{xx}|^2)dx\\
    &\leq Cs^{-1}\left(\|P_pv\|^2+\int_{\mathbb{T}}(s^8|v|^2+s^4|v_x|^2+s^4|v_{xx}|^2)\right)\\
    &= C\left(s^{-1}\|P_pv\|^2+\int_{\mathbb{T}}(s^7|v|^2+s^3|v_x|^2+s^3|v_{xx}|^2)\right),
\end{align*}
for some $C>0$. The above estimate and \eqref{precarleelliti} leads to
\begin{equation}\label{quasielitica}
\begin{array}{l}
\displaystyle\int_{\mathbb{T}}\{s^{-1}|v_{xxxx}|^2+s|v_{xxx}|^2+s^3|v_{xx}|^2+s^5|v_x|^2+s^7|v|^2\}dx\\
\leq C\Big(s^{-1}\|P_pv\|^2+\displaystyle\int_{\omega_0}(s^7|v|^2+s^3|v_x|^2+s^3|v_{xx}|^2)\\
+\displaystyle\int_{\omega_0}(s|v_{xxx}|^2+s^3|v_{xx}|^2+s^5|v_x|^2+s^7|v|^2)\Big)\\
\leq C\left(s^{-1}\|P_pv\|^2+\displaystyle\int_{\omega_0}(s|v_{xxx}|^2+s^3|v_{xx}|^2+s^5|v_x|^2+s^7|v|^2)\right),
\end{array}
\end{equation}
where $C$ does not depend on $s$ and $v$. In order to absorb the terms $v_x$ and $v_{xxx}$ above, let us consider $\xi\in C_{0}^{\infty}(\omega)$
with $0\leq \xi\leq 1$ for $x\in\omega_0$. Then, observe that
\begin{align*}
    \int_{\omega_0}|v_x|^2dx&\leq \int_{\omega}\xi|v_x|^2dx=\int_{\omega}\xi v_xv_xdx
    =-\int_{\omega}(\xi_xv_x+\xi v_{xx})vdx\\
    &=\frac{1}{2}\int_{\omega}\xi_{xx}v^2dx-\int_{\omega}\xi v_{xx}vdx,
\end{align*}
which allows us to deduce that
\begin{align}\label{condideriva1}
    2\int_{\omega_0}s|v_x|^2dx\leq \|\xi_{xx}\|_{L^{\infty}(\mathbb{T})}\int_{\omega}s|v|^2dx+\kappa \int_{\omega}s^{-1}|v_{xx}|^2dx+\int_{\omega}s^3|v|^2dx
\end{align}
and
\begin{align*}
     \int_{\omega_0}|v_{xxx}|^2dx&\leq \int_{\omega}\xi|v_{xxx}|^2dx=\int_{\omega}\xi v_{xxx}v_{xxx}dx
    =-\int_{\omega}(\xi_xv_{xxx}+\xi v_{xxxx})v_{xx}dx\\
    &=\frac{1}{2}\int_{\omega}\xi_{xx}v_{xx}^{2}dx-\int_{\omega}\xi v_{xxxx}v_{xx}dx.
\end{align*}
Consequently,
\begin{align}\label{condideriva2}
    2\int_{\omega_0}s|v_{xxx}|^2dx\leq \|\xi_{xx}\|_{L^{\infty}(\mathbb{T})}\int_{\omega}s|v_{xx}|^2dx+\kappa\int_{\omega}s^{-1}|v_{xxxx}|^2dx+\kappa^{-1}\int_{\omega}s^3|v_{xx}|^2dx,
\end{align}
for any $\kappa>0$. Combining \eqref{quasielitica}, \eqref{condideriva1} and \eqref{condideriva2} with $\kappa$ small enough and $s\geq s_0$ sufficiently large, for some constant $C$ that does not depend
on $s$ and $v$ we get
\begin{align}\label{Carlemprep}
    \int_{\omega}\{s^{-1}|v_{xxxx}|^2+s|v_{xxx}|^2+s^3|v_{xx}|^2+s^5|v_x|^2&+s^7|v|^2\}dx\leq \\ \nonumber
    &C\left( \|e^{s\psi}Pu\|^2+\int_{\omega}(s^7|v|^2+s^3|v_{xx}|^2)dx\right).
\end{align}
Replacing $v$ by $e^{s\psi}u$ in \eqref{Carlemprep} we obtain \eqref{carlemanelitica}. The proof of Lemma \ref{Carleli} is complete.
\end{proof}

We can now complete the proof of Proposition \ref{carleman}. Let $u\in L^2(0,T;H^4(\mathbb{T}))$ satisfying \eqref{carleq} and let $w=u-b_1u_{xx}+bu_{xxxx}\in L^2(0,T;L^2(\mathbb{T}))$.
Then,
$$w_t+\frac{a}{b}w_x=(\frac{a}{b}-q)u_x-(\frac{ab_1}{b}+p)u_{xxx}-ru_{xx}\in L^2(0,T;L^2(\mathbb{T})).$$
Combining \eqref{elitica}, \eqref{transpor} (multiplied by $e^{-2s\rho c^2t^2}$ and next integrated over $(0,T)$) and  \eqref{carlemanelitica}, it follows that, for $s\geq s_1$, the
following estimate holds
\begin{equation}
\begin{array}{l}
    \displaystyle\int_{0}^{T}\int_{\mathbb{T}}[s|u_{xxx}|^2+s^3|u_{xx}|^2+s^5|u_x|^2+s^7|u|^2+s|u-b_1u_{xx}+bu_{xxxx}|^2]e^{2s\varphi}dxdt\\ \nonumber
    +\displaystyle\int_{\mathbb{T}}[s|u-b_1u_{xx}+bu_{xxxx}|^2e^{2s\varphi}]_{t=0}dx \\ \nonumber
    \leq C\displaystyle\int_{0}^{T}\int_{\mathbb{T}}[|u_{xxxx}|^2+|(\frac{a}{b}-q)u_x-(\frac{ab_1}{b}+p)u_{xxx}-ru_{xx}|^2]e^{2s\varphi}dxdt\\ \nonumber
    +C\displaystyle\int_{0}^{T}\int_{\omega}[s|u-b_1u_{xx}+bu_{xxxx}|^2+s^7|u|^2+s^3|u_{xx}|^2]]e^{2s\varphi}dxdt.\nonumber
\end{array}
\end{equation}
Then, choosing $s_2\geq s_1$ and $C_2>C$ large enough, we obtain \eqref{carlecomple} for any $s\geq s_2$ and any $u\in L^2(0,T;H^4(\mathbb{T}))$ satisfying \eqref{carleq}.
\end{proof}

Finally, we can prove Theorem \ref{ucp}.

\begin{itemize}
\item If $u\in L^2(0,T;H^4(\mathbb{T}))$ satisfies \eqref{carleq} and \eqref{contuni}, it follows from \eqref{carlecomple} that $u=0$ in $\mathbb{T}\times(0,T)$.

\item If $u\in L^{\infty}(0,T;H^3(\mathbb{T}))$, then $u$ and $w=u-b_1u_{xx}+bu_{xxxx}$ are not regular enough to apply Lemmas \ref{Carleli} and \ref{carltra}. Therefore,
we first smooth them by using the following result (see, for instance, Proposition 1.4.29 in \cite{cazenave} and \cite{ro-zhan}).
\end{itemize}

\begin{proposition} If $X$ is a Banach space and $g\in L^p(0,T,X)$, with $1\leq p\leq\infty$, then for any $h > 0$ the function given by
$$ g^{[h]}(x,t)=\frac{1}{h}\int_{t}^{t+h}g(x,s)ds,$$
satisfies
\begin{itemize}
\item [{(i)}] $g^{[h]}\in W^{1,p}(0,T-h;X)$,
\item [{(ii)}] $\|g^{[h]}\|_{L^p(0,T-h;X)}\leq\|g\|_{L^p(0,T;X)}$,
\item [{(iii)}] $g^{[h]}\rightarrow g  \text{ in } L^p(0,T';X), \text{ as } h\rightarrow 0,$ for $p<\infty$ and $T'<T$.
\end{itemize}
\end{proposition}

Under the above conditions, pick any $T'\in (\frac{2b\pi}{a},T)$. Then, for any positive number $h<h_0=T-T'$, the function $u^{[h]}\in W^{1,\infty}(0,T',H^3(\mathbb{T}))$ and
solves the equation
\begin{equation}\label{equaaprox}
\begin{array}{l}
\vspace{2mm} u_{t}^{[h]}-b_1u_{txx}^{[h]}+bu_{txxxx}^{[h]}+au_{xxxxx}^{[h]}\\
    \qquad\qquad+(q(u)u_x)^{[h]}+(p(u)u_{xxx})^{[h]}+(r(u)u_{xx})^{[h]} =0\ \text{in}\ L^{\infty}(0,T';H^{-2}(\mathbb{T})),
    \end{array}
    \end{equation}
where $v_{t}^{[h]}$ denote $(v^{[h]})_t$, $v_{x}^{[h]}$ denote $(v^{[h]})_x$, etc. Moreover,
\begin{equation}\label{UCPapro}
\begin{array}{l}
    u^{[h]}(x,t)=0\ \ \ \text{in} \ \ \ (x,t)\in\omega\times(0,T').
\end{array}
\end{equation}
From \eqref{equaaprox}, we infer that
\begin{align*}
 &u_{xxxxx}^{[h]}\\
 &=a^{-1}(-u_{t}^{[h]}+b_1u_{txx}^{[h]}-bu_{txxxx}^{[h]}-(q(u)u_x)^{[h]}-(p(u)u_{xxx})^{[h]}-(r(u)u_{xx})^{[h]}).
\end{align*}
Hence, since the right hand side of the above identity belongs to $L^{\infty}(0,T';H^{-1}(\mathbb{T}))$,
\begin{align}\label{aproxregu}
    u^{[h]}\in L^{\infty}(0,T';H^4(\mathbb{T})).
\end{align}
This yields, with \eqref{elitica} and \eqref{transpor},
\begin{align}
    &w^{[h]}=u^{[h]}-b_1u_{xx}^{[h]}+bu_{xxxx}^{[h]}\in L^{\infty}(0,T';L^2(\mathbb{T})) \label{elitiapro}\\
    &w_{t}^{[h]}+\frac{a}{b}w_{x}^{[h]}=[(\frac{a}{b}-q)u_x]^{[h]}-[(\frac{ab_1}{b}+p)u_{xxx}]^{[h]}-[ru_{xx}]^{[h]}\in L^{\infty}(0,T';L^2(\mathbb{T}))\label{transapro}.
\end{align}
Thus, from \eqref{equaaprox}-\eqref{transapro}, Lemma \ref{Carleli} and \ref{carltra}, we obtain constants $s_1>0$ and $C_1>0$, such that, for all $s\geq s_1$ and all $h\in(0,h_0)$, we have
\begin{equation}\label{aproxcarl}
\begin{array}{l}
    \displaystyle\int_{0}^{T'}\int_{\mathbb{T}}[s|u_{xxx}^{[h]}|^2+s^3|u_{xx}^{[h]}|^2+s^5|u_{x}^{[h]}|^2+s^7|u^{[h]}|^2+s|u_{xxxx}^{[h]}|^2]e^{2s\varphi}dxdt\\
    \leq C\displaystyle\int_{0}^{T'}\int_{\mathbb{T}}[|u_{xxxx}^{[h]}|^2+\left|[(\frac{a}{b}-q)u_x]^{[h]}-[(\frac{ab_1}{b}+p)u_{xxx}]^{[h]}-[ru_{xx}]^{[h]}\right|^2]e^{2s\varphi}dxdt\\
    \leq C\displaystyle\int_{0}^{T'}\int_{\mathbb{T}}[|u_{xxxx}^{[h]}|^2+\left|[(\frac{a}{b}-q)u_x]^{[h]}\right|^2+\left|[(\frac{ab_1}{b}+p)u_{xxx}]^{[h]}\right|^2+\left|[ru_{xx}]^{[h]}\right|^2]e^{2s\varphi}dxdt\\
    \leq C\displaystyle\int_{0}^{T'}\int_{\mathbb{T}}[|u_{xxxx}^{[h]}|^2+\left|[(\frac{a}{b}-q)u_{x}^{[h]}\right|^2+\left|[(\frac{ab_1}{b}+p)u_{xxx}^{[h]}\right|^2+\left|[ru_{xx}^{[h]}\right|^2]e^{2s\varphi}dxdt\\
    +\displaystyle\int_{0}^{T'}\int_{\mathbb{T}}\left|[(\frac{a}{b}-q)u_x]^{[h]}-(\frac{a}{b}-q)u_{x}^{[h]}\right|^2e^{2s\varphi}dxdt\\
    +\displaystyle\int_{0}^{T'}\int_{\mathbb{T}}\left|[(\frac{ab_1}{b}+p)u_{xxx}]^{[h]}-(\frac{ab_1}{b}+p)u_{xxx}^{[h]}\right|^2e^{2s\varphi}dxdt\\
    +\displaystyle\int_{0}^{T'}\int_{\mathbb{T}}\left|[ru_{xx}]^{[h]}-ru_{xx}^{[h]}\right|^2e^{2s\varphi}dxdt.
\end{array}
\end{equation}
Comparing the powers of $s$ in \eqref{aproxcarl}, we obtain the following estimate for $s\geq s_3>s_1,h\in(0,h_0)$ and some constant $C_3>C_1$ (that does depend of $s,h$):
\begin{equation}\label{aproxcarl-1}
\begin{array}{l}
    \displaystyle\int_{0}^{T'}\int_{\mathbb{T}}[s|u_{xxx}^{[h]}|^2+s^3|u_{xx}^{[h]}|^2+s^5|u_{x}^{[h]}|^2+s^7|u^{[h]}|^2+s|u_{xxxx}^{[h]}|^2]e^{2s\varphi}dxdt\\
    \leq C\displaystyle\int_{0}^{T'}\int_{\mathbb{T}}|[(\frac{a}{b}-q)u_x]^{[h]}-(\frac{a}{b}-q)u_{x}^{[h]}|^2e^{2s\varphi}dxdt\\
    +C\displaystyle\int_{0}^{T'}\int_{\mathbb{T}}|[(\frac{ab_1}{b}+p)u_{xxx}]^{[h]}-(\frac{ab_1}{b}+p)u_{xxx}^{[h]}|^2e^{2s\varphi}dxdt\\
    +C\displaystyle\int_{0}^{T'}\int_{\mathbb{T}}|[ru_{xx}]^{[h]}-ru_{xx}^{[h]}|^2e^{2s\varphi}dxdt.
\end{array}
\end{equation}
In order to pass the right hand side of \eqref{aproxcarl-1} to the limit, we observe that, as $h\rightarrow0$,
\begin{align*}
  &[(\frac{a}{b}-q)u_x]^{[h]}\rightarrow (\frac{a}{b}-q)u_x \ \ \ \text{in} \ \ \ L^2(0,T';L^2(\mathbb{T})),\\
  &(\frac{a}{b}-q)u_{x}^{[h]}\rightarrow (\frac{a}{b}-q)u_x \ \ \ \text{in} \ \ \ L^2(0,T';L^2(\mathbb{T})),\\
  &[(\frac{ab_1}{b}+p)u_{xxx}]^{[h]}\rightarrow (\frac{ab_1}{b}+p)u_{xxx}\ \ \ \text{in} \ \ \ L^2(0,T';L^2(\mathbb{T})),\\
  &(\frac{ab_1}{b}+p)u_{xxx}^{[h]}\rightarrow (\frac{ab_1}{b}+p)u_{xxx}\ \ \ \text{in} \ \ \ L^2(0,T';L^2(\mathbb{T})),\\
  &[ru_{xx}]^{[h]}\rightarrow ru_{xx}\ \ \ \text{in} \ \ \ L^2(0,T';L^2(\mathbb{T})),\\
  &ru_{xx}^{[h]}\rightarrow ru_{xx}\ \ \ \text{in} \ \ \ L^2(0,T';L^2(\mathbb{T})),
\end{align*}
while $e^{2s_3\varphi}\in L^{\infty}(\mathbb{T}\times(0,T'))$. Consequently, for fixed $s$, we get
\begin{align*}
 &\int_{0}^{T'}\int_{\mathbb{T}}|[(\frac{a}{b}-q)u_x]^{[h]}-(\frac{a}{b}-q)u_{x}^{[h]}|^2e^{2s_3\varphi}dxdt\rightarrow0, \ \text{as}\ h\rightarrow 0,\\
 &\int_{0}^{T'}\int_{\mathbb{T}}|[(\frac{ab_1}{b}+p)u_{xxx}]^{[h]}-(\frac{ab_1}{b}+p)u_{xxx}^{[h]}|^2e^{2s_3\varphi}dxdt\rightarrow0, \ \text{as}\ h\rightarrow 0,\\
 & \int_{0}^{T'}\int_{\mathbb{T}}|[ru_{xx}]^{[h]}-ru_{xx}^{[h]}|^2e^{2s_3\varphi}dxdt\rightarrow0, \ \text{as}\ h\rightarrow 0.
\end{align*}
The convergences above and \eqref{aproxcarl-1} allow us to conclude that, as $h\rightarrow 0$,
\begin{align*}
    \int_{0}^{T'}\int_{\mathbb{T}}|u^{[h]}|^2e^{2s_3\varphi}dxdt\rightarrow 0.
\end{align*}
Since $u^{[h]}\rightarrow u$ in $L^2(0,T';L^2(\mathbb{T}))$,
\begin{align*}
    \int_{0}^{T'}\int_{\mathbb{T}}|u^{[h]}|^2e^{2s_3\varphi}dxdt\rightarrow\int_{0}^{T'}\int_{\mathbb{T}}|u|^2e^{2s_3\varphi}dxdt,
\end{align*}
as $h\rightarrow0$. Therefore, we can conclude that $u=0$ in $\mathbb{T}\times(0,T')$. As $T'$ may be taken arbitrarily close to $T$, we infer that $u=0$ in $\mathbb{T}\times(0,T)$.
\end{proof}
\section{Exponential stabilization}
Having the unique continuation result in hands, we derive the exponential decay of the solutions of \eqref{nlinitiesta-1} in the energy space $H^2(\mathbb{T})$, as $t\rightarrow\infty$.
This is done under suitable assumptions on the initial data.

Before going into the stabilization problem, let us explain the expression of the damping $\sigma$. We first write the linearized system of \eqref{nlinitiesta-1} as
\begin{equation}
\begin{cases}\label{cauchyest}
        u_t=Au+Bk,\\
    u(0)=u_0,
\end{cases}
\end{equation}
where $A=-(I-b_1\partial_{x}^{2}+b\partial_{x}^{4})^{-1}(\partial_x+a_1\partial_{x}^{3}+a\partial_{x}^{5})$, $k(t)=(I-b_1\partial_{x}^{2}+b\partial_{x}^{4})^{-1}h(t)\in L^2(0,T;H^s(\mathbb{T}))$ is a control input and
$$B=(I-b_1\partial_{x}^{2}+b\partial_{x}^{4})^{-1}\sigma(I-b_1\partial_{x}^{2}+b\partial_{x}^{4}).$$
We know from \cite{bau-pazo} that  $A$ is skew-adjoint in $H^s(\mathbb{T})$, and that \eqref{cauchyest} is exactly controllable in $H^s(\mathbb{T})$. Moreover, if we choose the simple feedback law (see \cite{LIU,SLEM})
\begin{align}\label{adjt}
k=-B^{*,s}u
\end{align}
the resulting closed-loop system
\begin{equation}\label{cauchyest1}
\begin{cases}
    u_t=Au-BB^{*,s}u,\\
    u(0)=u_0
\end{cases}
\end{equation}
is exponentially stable in $H^s(\mathbb{T})$, where $B^{*,s}$ denotes the adjoint of $B$ in $\mathcal{L}(H^s(\mathbb{T}))$. It is shown in the
Appendix that $B^{*,s}$ is given by
\begin{align}\label{adjunto B}
    B^{*,s}=(1-b\partial_{x}^2+b_{1}\partial_{x}^4)^{1-\frac{s}{2}}\sigma(1-b\partial_{x}^2+b_{1}\partial_{x}^4)^{\frac{s}{2}-1}.
\end{align}
We also deduce that
$$B^{*,2}u=\sigma u.$$

In order to make more precise the results stated above, let $\Tilde{A}=A-BB^{*,2},$ where $(BB^{*,2})u=(I-b_1\partial_{x}^{2}+b\partial_{x}^{4})[\sigma(I-b_1\partial_{x}^{2}+b\partial_{x}^{4})(\sigma u)]$. Since $BB^{*,2}\in \mathcal{L}(H^s(\mathbb{T}))$ and $A$ is skew adjoint in $H^s(\mathbb{T})$, $\Tilde{A}$ is the infinitesimal generator of a group $\{W_a(t)\}_{t\in\mathbb{R}}$ on $H^s(\mathbb{T})$ (See \cite[Theorem 3.4]{phill}).

Then, we have the well-posedness and the following exponentially stabilization result for \eqref{cauchyest}.

\begin{lemma}\label{estalinear}
Let $\sigma\in C^{\infty}(\mathbb{T})$ with $\sigma\neq0$. Then, there exist a constant $\beta>0$, such that, for $s\geq 2$, one can find a constant $C_s>0$ for which the following holds for all $u_0\in H^s(\mathbb{T})$:
\begin{align}\label{linestab}
    \|W_a(t)u_0\|_{H^s}\leq C_se^{-\beta t}\|u_0\|_{H^s},\ \ \forall\ \ \ t\geq0.
\end{align}
\end{lemma}
Adding the feedback law $k=-B^{*,2}u=-\sigma u$ in the nonlinear equation \eqref{Kdv-BBM-1} gives the closed-loop system
\eqref{nlinitiesta-1}. Then, Lemma \ref{estalinear} and the analysis developed in \cite{pas-pan} give the following well-posedness result
in the space $H^s(\mathbb{T})$ for $s\geq1$:
\begin{theorem}
Let $s\geq 2$ and $T>0$ be given. For any $u_0\in H^s(\mathbb{T})$, the system \eqref{nlinitiesta-1} admits a unique solution $u\in C([0,T];H^s(\mathbb{T}))$.
\end{theorem}

The next steps are devoted to show that \eqref{nlinitiesta-1} is globally exponentially stable in the space $H^2(\mathbb{T})$. In order to do that, were start by proving
the following observability inequality.

\begin{proposition}\label{estiestab}
Let $R_0>0$ and $r>0$ be given. Then, there exist two positive number $T$ and $\theta$, such that, for any $u_0\in H^3(\mathbb{T})$ satisfying
\begin{align*}
    0<r\leq \|u_0\|_{H^2(\mathbb{T})}\, \mbox{ and }\,\|u_0\|_{H^3(\mathbb{T})}\leq R_0,
\end{align*}
the corresponding solution $u$ of \eqref{nlinitiesta-1} satisfies
\begin{align}\label{observinequ}
    \|u_0\|_{H^2(\mathbb{T})}^{2}\leq \theta\int_{0}^{T}\|\sigma u(t)\|_{H^2(\mathbb{T})}^{2}dt.
\end{align}
\end{proposition}
\begin{proof}
Let $T>\frac{2b\pi}{|a|}$. We argue by contradiction and suppose that \eqref{observinequ} is not true. In this case, for any $n\geq 1$, \eqref{nlinitiesta-1}
admits a solution $u_n\in C([0,T];H^2(\mathbb{T}))$ satisfying
\begin{align}\label{limibol1}
     0<r\leq \|u_{0,n}\|_{H^2(\mathbb{T})}\, \mbox{ and }\,\|u_{0,n}\|_{H^3(\mathbb{T})}\leq R_0
\end{align}
and
\begin{align}\label{observinequ1}
   \int_{0}^{T}\|\sigma u_n(t)\|_{H^2(\mathbb{T})}^{2}dt\leq\frac{1}{n} \|u_{0,n}\|_{H^2(\mathbb{T})}^{2},
\end{align}
where $u_{0,n}=u_n(0)$. Since $\alpha_n=\|u_{0,n}\|_{H^2(\mathbb{T})}\leq \|u_{0,n}\|_{H^3(\mathbb{T})} \leq R_0$, we can choose a subsequence of $(\alpha_n)$, still denoted by $(\alpha_n)_{n\in\mathbb{N}}$, such that $\lim_{n\rightarrow\infty}\alpha_n=\alpha$. Observe that, from \eqref{observinequ1}, we obtain $\alpha_n>0$.

Following the notation introduced above, we introduce the function $v_n=\displaystyle\frac{u_n}{\alpha_n}$, for all $n\geq 1$. Then, $v_n$ satisfies
\begin{equation}\label{novaequa}
\begin{array}{l}
\vspace{2mm}v_{n,t}+v_{n,x}-b_1v_{n,txx}+a_1v_{n,xxx}+bv_{n,txxxx}+av_{n,xxxxx}\\\nonumber
    \qquad\qquad+\displaystyle\frac{3}{2}\alpha_nv_{n}v_{n,x}-\gamma\alpha_n(v_{n}^{2})_{xxx}-\frac{7}{48}\alpha_n(v_{n,x}^{2})_x
    -\frac{1}{8}\alpha_{n}^{2}(v_{n}^{3})_x=-\sigma(I-b_1\partial_{x}^{2}+b\partial_{x}^{4})\sigma v_{n}\nonumber
\end{array}
\end{equation}
and
\begin{align}\label{sequezero}
    \int_{0}^{T}\|\sigma v_n(t)\|_{H^2(\mathbb{T})}^{2}dt\leq\frac{1}{n}.
\end{align}
Moreover,
\begin{align}\label{norma1}
    \|v_n(0)\|_{H^2(\mathbb{T})}=1.
\end{align}
Since
$$\|v_n(0)\|_{H^3(\mathbb{T})}=\left\|\frac{u_n(0)}{\alpha_n}\right\|_{H^3(\mathbb{T})}=\frac{\|u_n(0)\|_{H^3(\mathbb{T})}}{\alpha_n}\leq \frac{R}{r},$$
for all $n\in\mathbb{N}$, the sequence $(v_n)_{n\in\mathbb{N}}$ is bounded in $L^{\infty}(0,T;H^3(\mathbb{T}))$ and $(v_{n,t})_{n\in\mathbb{N}}$ is bounded in $L^{\infty}(0,T;H^2(\mathbb{T}))$. Then, from Aubin-Lions lemma, we deduce that $(v_n)_{n\in\mathbb{N}}$ is bounded in $C(0,T;H^s(\mathbb{T}))$ for $2<s<3$. Therefore, we can extract a subsequence of $(v_n)_{n\in\mathbb{N}}$, still denoted by $(v_n)_{n\in\mathbb{N}}$, such that
\begin{align}\label{convergen}
    &v_n\rightarrow v \ \ \ \  \text{in} \ \ \ \ C(0,T;H^s(\mathbb{T})),\\
      &v_n\rightarrow v \ \ \ \  \text{in} \ \ \ \ L^{\infty}(0,T;H^3(\mathbb{T}))\ \ \text{weak-*,\label{convergef}}
\end{align}
for some $v\in L^{\infty}(0,T;H^3(\mathbb{T}))\cap C([0,T];H^s(\mathbb{T}))$, for all $2<s<3$. Consequently, from \eqref{convergen}
and \eqref{convergef}, we have that
\begin{align*}
    &\alpha_nv_nv_{n,x}\rightarrow\alpha vv_x\ \ \text{in}\ \  L^{\infty}(0,T,L^2(\mathbb{T}))\ \ \text{weak-*},\\
    &\alpha_n(v_{n}^{2})_{xxx}\rightarrow\alpha (v^2)_{xxx}\ \ \text{in}\ \  L^{\infty}(0,T,L^2(\mathbb{T}))\ \ \text{weak-*},\\
    &\alpha_n(v_{n,x}^{2})_x\rightarrow\alpha (v_{x}^{2})_x\ \ \text{in}\ \  L^{\infty}(0,T,L^2(\mathbb{T}))\ \ \text{weak-*},\\
    &\alpha_{n}^{2}(v_{n}^{3})_x\rightarrow\alpha^2 (v^3)_x\ \ \text{in}\ \  L^{\infty}(0,T,L^2(\mathbb{T}))\ \ \text{weak-*}.
\end{align*}
Furthermore, by \eqref{sequezero},
\begin{align}
    \int_{0}^{T}\|\sigma v\|_{H^2(\mathbb{T})}^{2}dt\leq \liminf_{n\rightarrow\infty}\int_{0}^{T}\|\sigma v_n\|_{H^2(\mathbb{T})}^{2}dt=0.
\end{align}
Thus, $v$ solves
\begin{equation}
\begin{array}{l}
\vspace{2mm} v_t+v_x-b_1v_{txx}+a_1v_{xxx}+bv_{txxxx}+av_{xxxxx}\\\nonumber
+\displaystyle\quad\frac{3}{2}\alpha vv_x-\gamma\alpha(v^2)_{xxx}-\frac{7}{48}(v_{x}^{2})_x-\frac{1}{8}\alpha^2(v^3)_x=-\sigma(I-b_1\partial_{x}^{2}+b\partial_{x}^{4})\sigma v, \ \ \ \ \ (x,t)\in\mathbb{T}\times(0,T),\nonumber
  \end{array}
\end{equation}
and, in addition,
$$v=0, \ \ \text{in}\ \ \ \ \omega\times(0,T).$$
From the UCP proved in Theorem \ref{ucp} we conclude that $v\equiv 0$ in $\mathbb{T}\times(0,T)$.

In order to obtain a contradiction, we first claim that $(v_n)_{n\in\mathbb{N}}$ is linearizable in the sense of \cite[Proposition 9]{deh-ger}. This is to say that, if $(w_n)_{n\in\mathbb{N}}$
denotes the sequence of solutions of the linear higher-order KdV-BBM equation with the same initial data as follows
\begin{equation}\label{equawn}
\begin{array}{l}
\vspace{2mm}w_{n,t}+w_{n,x}-b_1w_{n.txx}+a_1w_{xxx}+bw_{txxxx}+aw_{xxxxx}=-\sigma(I-b_1\partial_{x}^{2}+b\partial_{x}^{4})[\sigma w_n],\\
    w_n(x,0)=v_n(x,0),
\end{array}
\end{equation}
then
\begin{align}\label{dn}
    \sup_{0\leq t\leq T}\|v_n(t)-w_n(t)\|_{H^2(\mathbb{T})}\rightarrow0,\ \ \ \ \text{as}\ \ \ n\rightarrow\infty.
\end{align}
Indeed, if $d_n=v_n-w_n$, then $d_n $ solves
\begin{align*}
    &d_{n,t}+d_{n,x}-b_1d_{n,txx}+a_1d_{n,xxx}+bd_{txxxx}+ad_{xxxxx}=\\
    &-\sigma(I-b_1\partial_{x}^{2}+b\partial_{x}^{4})[\sigma d_n]-\frac{3}{2}\alpha_nv_nv_{n,x}-\gamma\alpha_n(v_{n}^{2})_{xxx}+\frac{7}{48}\alpha_n(v_{n,x}^{2})_x+\frac{1}{8}\alpha_{n}^{2}(v_{n}^{3})_x\\
    &d_n(0)=0.
\end{align*}
Since $\|W_a(t)\|_{\mathcal{L}(H^2(\mathbb{T}))}\leq Me^{\omega t}\leq Me^{\omega T}$, with $\omega,M>0$, from Duhamel formula we have that, for $t\in[0,T]$,
\begin{align*}
   \|d_n(t)\|_{H^2(\mathbb{T})}&\leq Me^{\omega T}\left( \int_{0}^{T}\|(I-b_1\partial_{x}^{2}+b\partial_{x}^{4})^{-1}\frac{3}{2}\alpha_nv_nv_{n,x}\|_{H^2(\mathbb{T})}dt\right.\\
   &+\int_{0}^{T}\|(I-b_1\partial_{x}^{2}+b\partial_{x}^{4})^{-1}\gamma\alpha_n(v_{n}^{2})_{xxx}\|_{H^2(\mathbb{T})}dt\\
   &+\int_{0}^{T}\|(I-b_1\partial_{x}^{2}+b\partial_{x}^{4})^{-1}\frac{7}{48}\alpha_n(v_{n,x}^{2})_x\|_{H^2(\mathbb{T})}dt\\
   &\left.+\int_{0}^{T}\|(I-b_1\partial_{x}^{2}+b\partial_{x}^{4})^{-1}\frac{1}{8}\alpha_n(v_{n}^{3})_x\|_{H^2(\mathbb{T})}dt\right).
\end{align*}
The above estimate combined with \eqref{convergen}-\eqref{convergef} and the fact that $v\equiv 0$ give us \eqref{dn}.

The next steps are devoted to prove that $\|v_n(0)\|_{H^2(\mathbb{T})}^{2}\rightarrow 0$, as $n\rightarrow\infty$. In fact,
from Lemma \ref{estalinear} we have that
\begin{align}\label{above}
    \|w_n(t)\|_{H^2(\mathbb{T})}\leq c_1e^{-\beta t} \|w_n(0)\|_{H^2(\mathbb{T})}, \ \ \ \text{for all}\ \ \ t\geq0,
\end{align}
and from  the energy identity for \eqref{equawn}, we get
\begin{align*}
     \|w_n(t)\|_{H^2(\mathbb{T})}^{2}- \|w_n(0)\|_{H^2(\mathbb{T})}^{2}=-2\int_{0}^{T} \|\sigma w_n(t)\|_{H^2(\mathbb{T})}^{2}dt
\end{align*}
or
\begin{align*}
    \|w_n(0)\|_{H^2(\mathbb{T})}^{2}= \|w_n(t)\|_{H^2(\mathbb{T})}^{2} + 2\int_{0}^{T} \|\sigma w_n(t)\|_{H^2(\mathbb{T})}^{2}dt.
\end{align*}
    Therefore, from \eqref{above} it follows that
    \begin{align}\label{final}
   \|w_n(0)\|_{H^2(\mathbb{T})}^{2}\leq 2(1-c_{1}^{2}e^{-2\beta T})^{-1}\left[\int_{0}^{T} \|\sigma w_n(t)-\sigma v_n(t)\|_{H^2(\mathbb{T})}^{2}dt+\int_{0}^{T} \|\sigma v_n(t)\|_{H^2(\mathbb{T})}^{2}dt\right].
\end{align}
Estimate \eqref{final} combined with \eqref{sequezero} and \eqref{dn} yields $\|v_n(0)\|_{H^2(\mathbb{T})}^{2}=\|w_n(0)\|_{H^2(\mathbb{T})}^{2}\rightarrow0$, as $n\rightarrow\infty$,
which contradicts \eqref{norma1}.
\end{proof}

The main result of this section reads as follows:

\begin{theorem}\label{estafinal}
Let $\sigma\in C^{\infty}(\mathbb{T})$ with $\sigma \neq0$, and $\beta>0$ be as given in Lemma \ref{estalinear}. Then, for any $R_0, r>0$, there exists a constant $C>0$, such that, for any $u_0\in H^3(\mathbb{T})$ with $0<r\leq \|u_0\|_{H^2(\mathbb{T})}\, \mbox{ and }\,\|u_0\|_{H^3(\mathbb{T})}\leq R_0$, the corresponding solution $u$ of \eqref{nlinitiesta-1} satisfies
\begin{align}
    \|u(\cdot,t)\|_{H^2(\mathbb{T})}\leq Ce^{-\beta t}\|u_0\|_{H^2(\mathbb{T})}.\nonumber
\end{align}
\end{theorem}
\begin{proof}
From Proposition \ref{estiestab} and the energy identity
$$\|u(t)\|_{H^2(\mathbb{T})}^{2}= \|u(0)\|_{H^2(\mathbb{T})}^{2}-2\int_{0}^{t} \|\sigma u(\tau)\|_{H^2(\mathbb{T})}^{2}d\tau,\ \ \ \ t\geq0,$$
we have
$$\|u(T)\|_{H^2(\mathbb{T})}^{2}\leq(1-2\theta^{-1}) \|u(0)\|_{H^2(\mathbb{T})}^{2}.$$
Thus,
$$\|u(kT)\|_{H^2(\mathbb{T})}^{2}\leq(1-2\theta^{-1})^k \|u(0)\|_{H^2(\mathbb{T})}^{2},\ \ \ k\in\mathbb{N},$$
which gives by the semigroup property
$$\|u(t)\|_{H^2(\mathbb{T})}^{2}\leq Ce^{-\kappa t} \|u(0)\|_{H^2(\mathbb{T})}^{2},\ \ \ \text{for all}\ \ \ t\geq0.$$
\end{proof}
\section{Appendix}\label{proof}

Proof of Lemma \ref{carltra}: (See \cite[Lemma 5.5]{ro-zhan}):

\vglue 0.2 cm

We first assume that $w\in H^1(\mathbb{T}\times(0,T))$. Let $v=e^{s\varphi}$ and $P=\partial_t+\frac{a}{b}\partial_x$. Then
\begin{align*}
    e^{s\varphi}Pw&=e^{s\varphi}P(e^{-s\varphi}v)\\
    &=(-s\varphi_tv-\frac{a}{b}s\varphi_xv)+(v_t+\frac{a}{b}v_x)\\
    &=P_nv+P_pv.
\end{align*}
It follows that
\begin{align}\label{normtran}
    \|e^{s\varphi}Pw\|_{L^2(\mathbb{T}\times(0,T))}^{2}= \|P_pv\|_{L^2(\mathbb{T}\times(0,T))}^{2}+ \|P_nv\|_{L^2(\mathbb{T}\times(0,T))}^{2}+2(P_pv,P_nv)_{L^2(\mathbb{T}\times(0,T))}^{2}.
\end{align}
After some integrations by parts in $t$ and $x$ in the last term in \eqref{normtran}, we obtain
\begin{align}\label{equavarp}
  2(P_pv,P_nv)_{L^2(\mathbb{T}\times(0,T))}^{2}=\int_{0}^{T}&\int_{\mathbb{T}}s(\varphi_{tt}+2\frac{a}{b}\varphi_{xt}+\frac{a^2}{b^2}\varphi_{xx})v^2dxdt\nonumber\\
  &-\int_{\mathbb{T}}s(\varphi_t+\frac{a}{b}\varphi_x)v^2|_{0}^{T}dx--\int_{0}^{T}\frac{a}{b}s(\varphi_t+\frac{a}{b}\varphi_x)v^2|_{0}^{2\pi}dt.
\end{align}
Using \eqref{condperi}-\eqref{varphi} and the fact that $v(0,t)=v(2\pi,t)$, we notice that the last term in \eqref{equavarp} is null. From \eqref{poli}-\eqref{varphi}, we infer that
\begin{align*}
    &\varphi_{tt}+2\frac{a}{b}\varphi{tx}+\frac{a^2}{b^2}\varphi_{xx}=2(1-\rho)\frac{a^2}{b^2}>0\ \ \ \text{for}\ \ (x,t)\in(\frac{\eta}{2},2\pi-\frac{\eta}{2})\times(0,T),\\
    &-(\varphi_t+\frac{a}{b}\varphi_x)\geq 2\frac{a}{b}(\frac{aT\rho}{b}-2\pi-\delta)>0\ \ \ \text{for}\ \ x\in(0,2\pi), t=T\\
    &\varphi_t+\frac{a}{b}\varphi_x\geq 2\frac{a}{b}\delta>0\ \ \ \text{for}\ \ x\in(0,2\pi), t=0.
\end{align*}
Thus,
\begin{align*}
    \int_{0}^{T}\int_{\mathbb{T}}s|v|^2dxdt+\int_{\mathbb{T}}s(|v|_{t=0}^{2}+|v_{t=T}^{2})dx\leq C\left(\int_{0}^{T}\int_{\mathbb{T}}|e^{s\varphi}Pw|^2dxdt+\int_{0}^{T}\int_{\omega}s|v|^2dxdt\right),
\end{align*}
which gives at once \eqref{Carlemantran} by replacing $v$ by $e^{s\varphi}w$. The proof of Lemma \ref{carltra} is achieved when $w\in H^1(\mathbb{T}\times(0,T))$.
We now claim that Lemma \ref{carltra} is still true when $w$ and $f$ are in $L^2(0,T;L^2(\mathbb{T}))$. Indeed, in the case $w\in C([0,T];L^2(\mathbb{T}))$, and if $(w_{n}^{0})$ and $f_n$ are two sequences in $H^1(\mathbb{T})$ and $L^2(0,T;H^1(\mathbb{T}))$, respectively, such that
\begin{align*}
    &w_{n}^{0}\rightarrow w(0)\ \ \ \text{in}\ \ L^2(\mathbb{T})\\
    &f^n\rightarrow f\ \ \ \text{in}\ \ L^2(0,T;L^2(\mathbb{T}),
\end{align*}
then the solution $w^n\in C([0,T];H^1(\mathbb{T}))$ of
\begin{align*}
    w_{t}^{n}+\frac{a}{b}w_{x}^{n}+f^n,\\
    w^n(0)=w_{0}^{n}
\end{align*}
satisfies $w^n\in H^1(\mathbb{T}\times(0,T))$ and $w^n\rightarrow w$ in $C([0,T];L^2(\mathbb{T}))$, so that we can apply \eqref{Carlemantran} to $w^n$ and next pass to the limit $n\rightarrow\infty$ in \eqref{Carlemantran}. The proof of Lemma \ref{carltra} is complete.

\vglue 0.2 cm

\noindent Proof of \eqref{adjunto B}:

\vglue 0.2 cm

Observe that 	
\begin{equation*}
c_{b,b_{1}}(1+bx^2+b_{1}x^4)^{\frac{s}{2}}\leq[(1+x^2)^{2}]^{\frac{s}{2}}\leq C_{b,b_{1}}(1+bx^2+b_{1}x^4)^{\frac{s}{2}},
\end{equation*}	
for $s\geq 2$ and some positive  constants  $c_{b,b_{1}},C_{b,b_{1}}.$ Then, we can define the following equivalent inner product in $H^{s}(\mathbb{T})$ as
\begin{equation*}
(u,v)_{s}=\int_{\mathbb{T}}(1+bx^2+b_{1}x^4)^{\frac{s}{2}}\mathcal{F}{u}(x)\overline{\mathcal{F}v}(x)d,
\end{equation*}	
where $\mathcal{F}\varphi$ denote the Fourier transform of $\varphi$. Hence, employing Plancherel Theorem, we get
\begin{equation*}
\begin{aligned}
(B\varphi,\psi)_{s}&= \int_{\mathbb{T}}(1+bx^{2}+b_{1}x^{4})^{\frac{s}{2}}\mathcal{F}[(1-b\partial_{x}^2+b_{1}\partial_{x}^4)^{-1}\sigma(x)(1-b\partial_{x}^2+b_{1}\partial_{x}^4)\varphi(x)] \ \overline{\mathcal{F}\psi}(x)dx
\\
&=\int_{\mathbb{T}}(1+bx^{2}+b_{1}x^{4})^{\frac{s}{2}-1}\mathcal{F}[\sigma(x)(1-b\partial_{x}^2+b_{1}\partial_{x}^4)\varphi(x)]  \overline{\mathcal{F}\psi}(x)dx
\\
&=\int_{\mathbb{T}}\mathcal{F}[\sigma(x)(1-b\partial_{x}^2+b_{1}\partial_{x}^4)\varphi(x)]  \overline{\mathcal{F}[(1-b\partial_{x}^2+b_{1}\partial_{x}^4)^{\frac{s}{2}-1}\psi(x)] }dx
\\
&=(\sigma(x)(1-b\partial_{x}^2+b_{1}\partial_{x}^4)\varphi,(1-b\partial_{x}^2+b_{1}\partial_{x}^4)^{\frac{s}{2}-1}\psi)_{L^2(\mathbb{T})}
\\
&=((1-b\partial_{x}^2+b_{1}\partial_{x}^4)\varphi,\sigma(x)(1-b\partial_{x}^2+b_{1}\partial_{x}^4)^{\frac{s}{2}-1}\psi)_{L^2(\mathbb{T})}
\\
&=\int_{\mathbb{T}}\mathcal{F}[(1-b\partial_{x}^2+b_{1}\partial_{x}^4)\varphi(x) ]\overline{\mathcal{F}[\sigma(x)(1-b\partial_{x}^2+b_{1}\partial_{x}^4)^{\frac{s}{2}-1}\psi(x)]}dx
\\
&=\int_{\mathbb{T}}(1+bx^2+b_{1}x^4)^{\frac{s}{2}}\mathcal{F}\varphi(x) (1+bx^2+b_{1}x^4)^{1-\frac{s}{2}}\overline{\mathcal{F}[\sigma(x)(1-b\partial_{x}^2+b_{1}\partial_{x}^4)^{\frac{s}{2}-1}\psi(x)]}dx
\\
&=\int_{\mathbb{T}}(1+bx^2+b_{1}x^4)^{\frac{s}{2}}\mathcal{F}\varphi(x) \overline{\mathcal{F}[(1-b\partial_{x}^2+b_{1}\partial_{x}^4)^{1-\frac{s}{2}}\sigma(x)(1-b\partial_{x}^2+b_{1}\partial_{x}^4)^{\frac{s}{2}-1}\psi(x)]}dx
\\
&=(\varphi,(1-b\partial_{x}^2+b_{1}\partial_{x}^4)^{1-\frac{s}{2}}\sigma(x)(1-b\partial_{x}^2+b_{1}\partial_{x}^4)^{\frac{s}{2}-1}\psi)_{s},
\end{aligned}
\end{equation*}
 for all $\varphi, \psi\in H^s(\mathbb{T})$. From the computation above we deduce \eqref{adjunto B}.


\begin{thebibliography}{99}


\bibitem{bau-pazo} G. J. Bautista and A. F. Pazoto,\, {\em A note on the control and stabilization of a higher-order water wave model},  Discrete Contin. Dyn. Syst. Ser. B, {\bf 28} (2023), 1513--1527.

\bibitem{bona-carvajal-panthee-scialom} J. L. Bona, X. Carvajal, M. Panthee and M. Scialom,\, {\em Higher-order Hamiltonian model for undirectional water waves,} J. Nonlinear
Sci. {\bf 28} (2018), 543--577.

\bibitem{BCS1} J. L. Bona, M.  Chen and J.-C. Saut,\,{\em Boussinesq equations and other systems for small-amplitude long waves in
			nonlinear dispersive media. I: Derivation and linear theory,} J. Nonlinear Sci., {\bf 12} (2002), 283--318.
		
		\bibitem{BCS2} J. L. Bona, M.  Chen and J.-C. Saut,\,{\em Boussinesq equations and other systems for small-amplitude long waves in
			nonlinear dispersive media. II: Nonlinear theory}, Nonlinearity, {\bf 17} (2004), 925--952.

\bibitem{pas-pan} X. Carvajal, M. Panthee and R. Pastr\'{a}n,\, {\em On the well-posedness, ill-posedness and norm-inflation for a higher order water wave model on a periodic domain,} Nonlinear Analysis. Theory, Methods and Applications.  {\bf 192} (2020), 111713, 22 pp.

\bibitem{cazenave} T. Cazenave and A. Haraux,\, {\em An Introduction to Semilinear Evolution Equation}, Oxford
Lecture Series in Mathematics and its Applications, {\bf 13}, The Clarendon Press, Oxford University Press, New York, 1998.

\bibitem{deh-ger} B. Dehman, P. Gérard and G. Lebeau,\, {\em Stabilization and control for the nonlinear Schrodinger equation on a compact surface,} Math. Z.  {\bf 254} (2006), 729--749.

\bibitem{LIU} K. Liu, {\em Locally distributed control and damping for the conservative systems}, SIAM J. Cont. Optim.,  {\bf 35} \, (1997), 1574--1590.

\bibitem{phill} R. S. Phillips\, {\em Perturbation theory for semi-groups of linear operators,} Transactions of the American Mathematical Society.  {\bf 74} (1953), 199--221.

\bibitem{rosier-zhang surv} L. Rosier and B.-Y. Zhang,\, {\em Control and stabilization of the Korteweg-de Vries equation: recent progresses}, J. Syst. Sci. Complex.
{\bf 22} (2009), 647-–682.

\bibitem{ro-zhan} L. Rosier and B.-Y. Zhang,\, {\em Unique continuation property and control for the Benjamin-Bona-Mahony equation on a periodic domain,}
J. Differential Equations {\bf 254} (2013), 141--178.

\bibitem{SLEM} M. Slemrod,\, {\em A note on complete controllability and stabilizability for linear control systems in a Hilbert space}, SIAM J. Cont. Optim {\bf 12}\,(1974), 500--508.

\end{thebibliography}
\end{document}